\documentclass{ifacconf}

\usepackage{amsmath,amssymb,amsfonts}
\usepackage{graphicx}
\usepackage{textcomp}
\usepackage{color}
\usepackage{wrapfig}
\usepackage{lipsum}
\usepackage{epstopdf}

\usepackage{amsopn}

\usepackage{enumerate}
\usepackage{cancel}
\usepackage{mathrsfs}
\usepackage{mathdots}
\usepackage{euscript}
\usepackage{amscd}
\usepackage{placeins}
\usepackage{natbib}

\usepackage{tikz}
\usetikzlibrary{arrows}
\usepackage{verbatim}
\usepackage{algorithm}
\usepackage{algpseudocode} 

\newcommand{\bmat}{\left[ \begin{matrix}}
	\newcommand{\emat}{\end{matrix} \right]}

\DeclareMathOperator*{\argmin}{argmin} 

\newcommand{\Rbb}{\mathbb R}

\newcommand{\Nbb}{\mathbb N}

\newcommand{\pb}{\mathbf  p}

\newcommand{\Ccal}{\mathcal{C}}
\newcommand{\Dcal}{\mathcal{D}}

\renewcommand{\alg}[1]{\begin{align} #1 \end{align}}
\newcommand{\nn}{\nonumber}

\newcommand\norm[1]{\lVert#1\rVert} 
\begin{document}
\begin{frontmatter}

\title{A kernel-based PEM estimator\\ for forward models} 

\author[Unipd]{Giulio Fattore}
\author[Unipd]{Marco Peruzzo}
\author[Unipd]{Giacomo Sartori}
\author[Unipd]{Mattia Zorzi} 
\address[Unipd]{Department of Information Engineering, University of Padova, Via Gradenigo 6/B, 35131 Padova, Italy\\\{\!fattoregiu,peruzzomar,zorzimat\}\!@dei.unipd.it, giacomo.sartori.4@studenti.unipd.it}

\begin{abstract}                
This paper addresses the problem of learning the impulse responses characterizing forward models by means of a regularized kernel-based Prediction Error Method (PEM). The common approach to accomplish that is to approximate the system with a high-order stable ARX model. However, such choice induces a certain undesired prior information in the system that we want to estimate. To overcome this issue, we propose a new kernel-based paradigm which is formulated directly in terms of the impulse responses of the forward model and leading to the identification of a high-order MAX
model. The most challenging step is the estimation of the kernel hyperparameters optimizing the marginal likelihood. The latter, indeed, does not admit a closed form expression. We propose a method for evaluating the marginal likelihood which makes possible the hyperparameters estimation. Finally, some numerical results showing the effectiveness of the method are presented. 
\end{abstract}

\begin{keyword}
System identification; Kernel-based PEM methods; Marginal likelihood optimization.
\end{keyword}

\end{frontmatter}

\section{Introduction}
System identification for linear-discrete time models is typically handled by using Prediction Error Methods (PEM) and parametric model classes \citep{LJUNG_SYS_ID_1999,SODERSTROM_STOICA_1988}. However, the choice of the model structure, which is usually accomplished by using Akaike Information Criterion (AIC) or Bayesian Information Criterion (BIC), see \cite{AKAIKE_1974,SCHWARZ_1978}, is the most challenging part.

An alternative approach to estimate a linear model resorts to regularized kernel-based PEM methods \citep{EST_TF_REVISITED_2012,PILLONETTO_DENICOLAO2010,doi:10.1080/00207179.2019.1578407,chen2013implementation,KHO1,KHO2,chen2016maximum,7495008,10384202}. In this context, we consider a model class containing high-order ARX models which are described through the predictor impulse responses. The search of the optimal model over this ``large'' model class is an ill-posed problem because only a finite set of measured data is available, however it can be made into a well-posed one introducing a penalty term that promotes models with specific features. Adopting the Bayesian viewpoint, the impulse responses are modeled as a zero-mean Gaussian random vector with an appropriate covariance (or kernel) matrix \citep{PILLONETTO_2011_PREDICTION_ERROR}. The latter characterizes the a-priori knowledge on the predictor impulse responses, such as the requirement for the impulse response to be absolutely integrable (i.e. the corresponding system is Bounded Input Bounded Output (BIBO) stable) and possess a certain level of smoothness. Many methods have been proposed in order to design kernels for system identification, see e.g.  
\cite{CHEN2018109,ZORZI2018125,marconato2017filter,zorzi2021second,chen2015kernel1,chen2015kernel2,FUJI1,FUJI2,ZORZI_MED}, or for network identification, see e.g. \cite{BSL,BKRON,BSL_CDC,CDC_KRON,KRON,REWEIGHTED,ALPAGO}. It is worth noting that the kernel depends on some hyperparameters which are usually inferred from the data by optimizing the so-called marginal likelihood \citep{RASMUSSEN_WILLIAMNS_2006}. 

\begin{figure}[!h]
    \centering
    \includegraphics[width=0.47\textwidth,height=12.5cm]{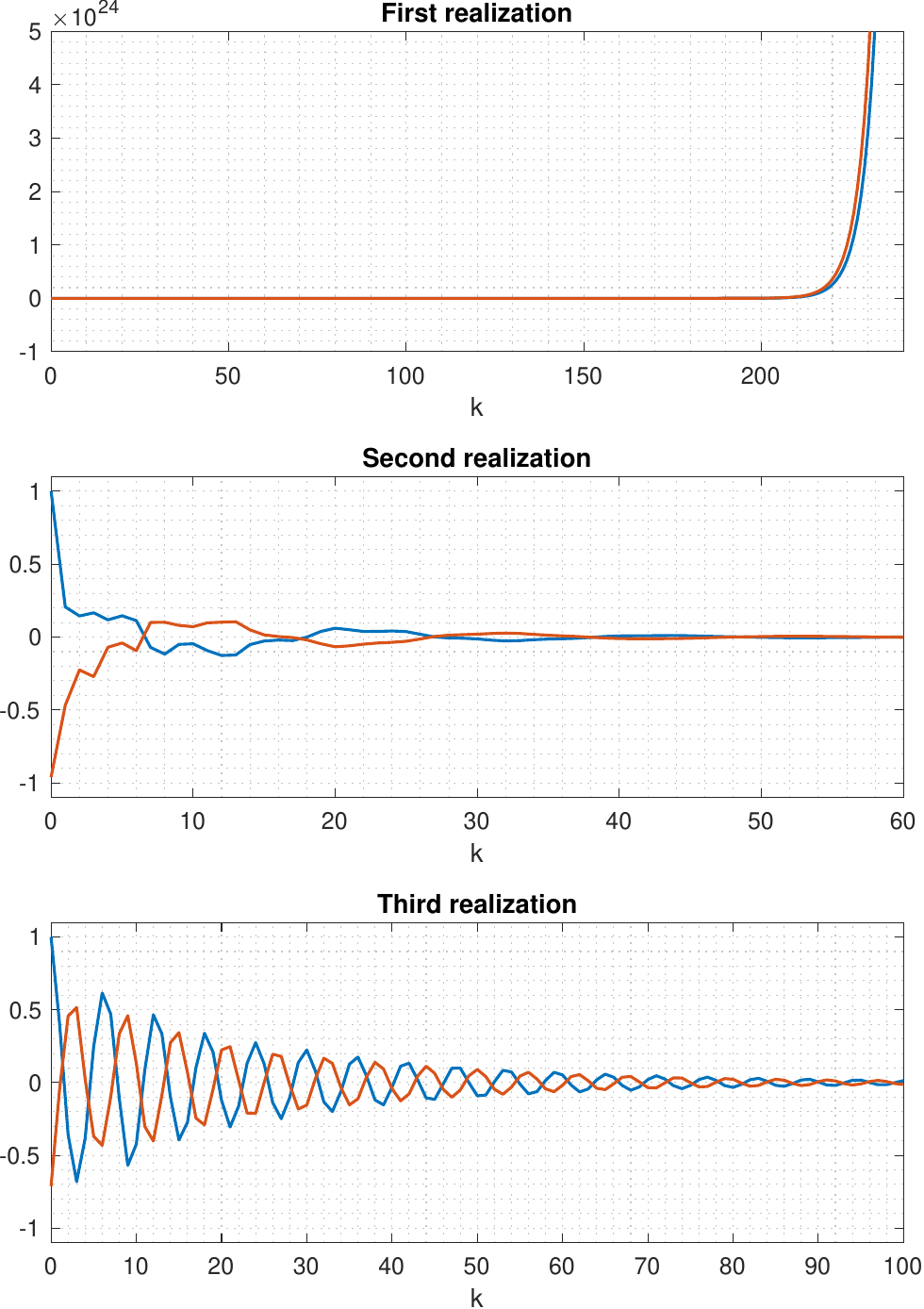}
    \caption{Three different realizations of the impulse responses $\{b_k,\; k\in\Nbb_0\}$ ({blue line}) and $\{c_k,\; k\in\Nbb_0\}$ ({orange line}) corresponding to $B(z)$ and ${C}(z)$.}
    \label{fig:correlation}
\end{figure}

In this paper we focus {on} the problem of learning the impulse responses characterizing the ``forward'' (or simulation) model: 
\begin{equation}\label{eq:nonpar_model}
 y(t)=B(z)u(t-1)+C(z)e(t)
\end{equation}
where $y,u,e$ denote the output, the input and the noise process, respectively; $B(z)$ and $C(z)$ are the (causal) transfer functions corresponding to the impulse responses of the forward model. In the literature, a well established method to solve the aforementioned problem is considering an high order ARX model 
approximation of the BIBO stable system (\ref{eq:nonpar_model}), see \cite[Section 2.2]{PILLONETTO2014}:
\begin{equation}\label{eq:arx}
y(t)=A(z)y(t)+F(z)u(t-1)+e(t)   
\end{equation}
where
\alg{A(z)\approx 1-\frac{1}{C(z)},\quad F(z)\approx \frac{B(z)}{C(z)} \label{eq:approx}}
are polynomial transfer functions and $A(z)$ is strictly causal.
Then, it is possible to use the kernel-based PEM method for identifying the model in (\ref{eq:arx}). However, it is not trivial to design the kernel for the predictor impulse responses corresponding to $A(z)$ and $F(z)$ in such a way to induce certain a priori information on the impulse responses of $B(z)$ and ${C}(z)$, i.e. the ones of the forward model. In addition, if we parametrize the impulse responses of $B(z)$ and $C(z)$ through the ones of $A(z)$ and $F(z)$, then the former impulse responses are Gaussian random vectors which are never independent. If such dependence does not exist, then such strategy adds a wrong a priori information in the system identification procedure. To better clarify that, we make an example.
We model the predictor impulse responses corresponding to $A(z)$ and ${F}(z)$ as two zero-mean Gaussian and independent processes with covariance function characterized by the Tuned Correlated (TC) kernel with decay rate equal to 0.6 and 0.5. We draw three realizations from these stochastic process{es}. Then, we compute the corresponding realizations of the forward impulse responses through the mapping in \eqref{eq:approx}, see Figure~\ref{fig:correlation}. 
In all the three cases the impulse responses share similar features and thus they are correlated: in the first realization, Figure~\ref{fig:correlation}(top panel), the forward impulse responses diverge; in the second realization, Figure~\ref{fig:correlation}(panel in the middle), the impulse responses are smooth; in the third realization, Figure~\ref{fig:correlation}(bottom panel), the impulse responses exhibit an oscillatory behavior. It is also worth noting that in the second and third cases the impulse responses have a similar decay rate.

To overcome the issues previously discussed, which arises using the ARX approximation,
we propose a new-kernel based method relying on a high-order MAX model structure. In this way we are able to directly formulate the problem in terms of impulse responses of the forward model, and thus design directly the corresponding priors. 
Within such setup, it is very difficult to find an analytic expression for the marginal likelihood, which is needed for estimating the kernel hyperparameters. As a consequence, the evaluation of the marginal is a non-trivial task since it requires the computation of an integral. Accordingly, we propose a method for evaluating the marginal likelihood and we show numerically its effectiveness.

The paper is organized as follows. In Section \ref{sec 2} the identification problem of forward systems as well as the corresponding Bayesian model are introduced. In Section \ref{Proposed Estimator section} we introduce the kernel-based PEM estimator; in particular, we show how to evaluate the marginal likelihood for estimating the hyperparameters of the kernel matrix. In Section \ref{sec num exp} some numerical experiment are presented in order to test the proposed estimator. Finally in Section \ref{conclusion} we draw the conclusions.

{\em Notation.} Given a matrix $K\in\mathbb R^{T\times T}$, $K^\top$ denotes the transpose of $K$, $(K)_{ij}$ denotes its entry in position $(i,j)$, while $|K|$ its determinant. If $K$ is positive definite, then the weighted norm of $\theta\in\mathbb R ^T$ corresponding to $K$ is defined as $\|\theta\|_{K^{-1}}=\sqrt{\theta^\top K^{-1}\theta}$. The notation $y\sim \cal N(\mu,K)$ means that $y$ is a Gaussian random vector with mean $\mu$ and covariance matrix $K$.

\section{Problem formulation}\label{sec 2}
Consider the nonparametric Single Input Single Output (SISO) linear and forward model 
\alg{\label{def_mod}y(t)=\sum_{k=0}^{\infty} b_k u(t-k-1)+\sum_{k=0}^\infty c_k e(t-k)}
where $c_0=1$; $y(t),u(t)$ denote the output and the input of the system at time $t$, respectively; $e(t)$ is white Gaussian noise with variance $\mathrm{var}(e(t))=\sigma^2$. The input can be deterministic or a stochastic process independent of $e$. We assume that the system is BIBO stable that is the impulse responses $\{ b_k,\;\ k\in\Nbb_0\}$, $\{ c_k,\;  k\in\Nbb_0\}$ are absolutely summable. Accordingly, we can choose $T\in\mathbb N$ sufficiently large in such a way that $|b_k|\approx 0$ and $|c_k|\approx 0$ for any $k>T$ and $T$ represents the practical length of such impulse responses. Therefore, the system can be approximated through the high-order MAX
model
\alg{\label{ARMAXhigh}y(t)=\sum_{k=0}^{T-1} b_k u(t-k-1)+\sum_{k=0}^T c_k e(t-k).}
Notice that the model above is completely characterized by $\sigma^2$ and the parameter vector
\begin{equation}
    \theta=\begin{bmatrix}
        b\\
        c
    \end{bmatrix}\in\Rbb^{2T}
\end{equation}
which is composed by the vectors 
\begin{equation} \label{eq:bc}
   b
   =\begin{bmatrix}
        b_0\\
        b_1\\
        \vdots\\
        b_{T-1}\\
    \end{bmatrix}, \
    c
    =\begin{bmatrix}
        c_1\\
        c_2\\
        \vdots\\
        c_T
    \end{bmatrix} .
\end{equation}

Assume to collect a dataset $\Dcal:=\{(y(t),u(t)),\; t=1\ldots N\}$ generated by 
(\ref{ARMAXhigh}). We want to estimate the impulse responses characterizing model (\ref{ARMAXhigh}) (that is indeed $\theta$) from $\Dcal$.

Adopting the Bayesian viewpoint, $b$ and $c$ are modeled as independent Gaussian random vectors with zero mean and covariance matrices (i.e. kernel matrices) $K_b\in \Rbb^{ T \times T}$ and $K_c\in \Rbb^{ T \times T}$, respectively. Therefore, we have $\theta \sim \mathcal{N}(0,K)$ where 
\begin{equation}\label{defK}
    K=\begin{bmatrix} 
       \lambda_b K_{b} & 0 \\ 0 &  \lambda_cK_{c} 
    \end{bmatrix}
\end{equation}
and $\lambda_b,\lambda_c\in(0,\infty)$ are the scale factors. Kernel matrices $K_b$ and $K_c$ are designed in such a way to induce some a priori knowledge about the impulse response, e.g. it should decay sufficiently fast to zero. In what follows we parametrize $K_b$ and $K_c$ through:

- the TC kernel \citep{EST_TF_REVISITED_2012}
\begin{equation}
 \label{defTC}   (K_{b})_{ij}=\beta_b^{max(i,j)}, \ (K_{c})_{ij}=\beta_c^{max(i,j)}
\end{equation}
where $\beta_b, \beta_c\in (0,1)$ are the decay rate hyperparameters;

- the second order Diagonal Correlated (DC2) kernel \citep{zorzi2021second}
\alg{\label{dDC2}
    (K_{b})_{ij}&=\beta_b^{max(i,j)}(1-(1-\beta_b)\alpha_b^{|i-j|+1})-\alpha_b^2\beta_b^{max(i,j)+1}, \nn\\ 
    (K_{c})_{ij}&=\beta_c^{max(i,j)}(1-(1-\beta_c)\alpha_c^{|i-j|+1})-\alpha_c^2\beta_c^{max(i,j)+1} }
where $\alpha_b, \alpha_c\in (0,1)$ are additional hyperparameters.

These kernels encode the a priori information that the impulse response decays to zero and has a certain degree of smoothness. More precisely, the DC2 kernel induces more smoothness on the impulse response than the TC kernel. Moreover, the smoothness degree in the DC2 kernels is tuned according to $\alpha_b$ and $\alpha_c$. If more smoothness is needed, then it is possible to use the Stable-Spline (SS) kernel \citep{PILLONETTO_DENICOLAO2010} or the second order TC (TC2) kernel \citep{zorzi2021second}.    

\section{Proposed Estimator }\label{Proposed Estimator section}
In this section we propose a kernel-based PEM estimator in order to estimate the parameter vector $\theta$ characterizing the high-order MAX 
model 
\eqref{ARMAXhigh} from the dataset $\Dcal$. First, we assume that $\sigma^2$ and the hyperparameters vector $\eta$ are given. For instance, $\eta=\{\lambda_b,\lambda_c,\beta_b,\beta_c,\alpha_b,\alpha_c\}$ in the case we use the DC2 kernel for $b$ and $c$. Then, in the next section we show how to estimate $\eta$ from  $\Dcal$.
In what follows all the probability densities are conditional on the input, but we omit this dependence to simplify the notation.

\subsection{MAP Estimator of $\theta$}
Let \alg{B(z)=\sum_{k=0}^{T-1} b_k z^{-k},\quad C(z)=1+\sum_{k=1}^T c_kz^{-k} }
denote the Z-transform of the truncated impulses $b$ and $c$. It is not difficult to see that the one-step ahead predictor of $y(t)$ under model (\ref{ARMAXhigh}) is 
\begin{equation} \label{hatytheta}
    \hat{y}_\theta(t|t-1)=[
    1-C(z)^{-1}]y(t)+\frac{B(z)}{C(z)}u(t-1) \nn
\end{equation}
where we made explicit its dependence on $\theta$. The average squared prediction error using $\Dcal$ is
\begin{equation} \label{Vn}
    V_N(\theta)=\frac{1}{N}\sum_{t=1}^N\varepsilon^2_\theta(t)
\end{equation}
where \begin{equation}
 \varepsilon_\theta(t)=y(t)-\hat{y}_\theta(t|t-1)
\end{equation}
denotes the one-step ahead prediction error at time $t$.
Then, the kernel-based PEM estimator is the solution of the optimization problem 
\begin{equation} \label{Theta_hat}
    \hat{\theta}_{\eta}= \argmin_{\theta \in \Theta}V_N(\theta)+ \frac{\sigma^2}{N}\norm{\theta}^2_{{K_\eta}^{-1}}  
\end{equation}
where: $\Theta$ is the set of the parameters $\theta$ for which the impulse response corresponding to $C(z)^{-1}$ is absolutely summable; $K_\eta$ is defined as in (\ref{defK}) where we made explicit the dependence on $\eta$. It is worth noting that the optimization problem (\ref{Theta_hat}) is not convex because the one-step ahead predictor defined in (\ref{hatytheta}) is not linear in $\theta$. More precisely, $\sigma^2N^{-1}\norm{\theta}^2_{{K_\eta}^{-1}}$ is a strictly convex function, while $V_N(\theta) $ does not. 
On the other hand, since the MAX structure is a particular case of ARMAX, as mentioned in \cite{LJUNG_SYS_ID_1999}, pag. 339: \lq\lq the global minimum (for $V_N(\theta)$) is usually found without too much problem for ARMAX models, see \cite{bohlin1971problem} for a detailed discussion\rq\rq.

Let $\mathrm y^-:=[\, y(0)\; y(-1) \,\ldots \,]^\top$ denote the (infinite length) past data vector of the output. We also define \alg{\mathrm y^+&:=[\, y(1)\; y(2) \,\ldots\, y(N)\,]^\top.\nn} Let $\pb(\mathrm y^+,\theta|\mathrm y^-,\eta,\sigma^2)$ denote joint probability of $\mathrm y^+$ and $\theta$ given $\mathrm y^-$, $\eta$ and $\sigma^2$. We 
assume that the past data $\mathrm y^-$ neither affects the a priori probability on $\theta$ nor carries information on $\eta$ and $\sigma^2$. Then, we have
\alg{\pb(\mathrm y^+,\mathrm y^-,\theta|\eta,\sigma^2)\approx \pb(\mathrm y^+|\theta,\mathrm y^-,\eta,\sigma^2)\pb(\theta|\eta,\sigma^2)\pb (\mathrm y^-). \nn}
Moreover, it is not difficult to see that
\alg{ \pb(\mathrm y^+|\theta,\mathrm y^-,\eta,\sigma^2) &=\prod_{t=1}^N\frac{1}{\sqrt{2\pi\sigma^2}}\exp\left(-\frac{1}{2\sigma^2} \varepsilon_\theta(t)^2 \right)\nn\\
\pb(\theta|\eta,\sigma^2)&=\frac{1}{\sqrt{(2\pi)^{2T}|{K_\eta}}|}\exp\left(-\frac{1}{2} \|\theta\|^2_{K_\eta^{-1}} \right).\nn}
Thus, the negative log-posterior is 
\alg{\label{neg_lik}
    \ell(\theta,\eta;\mathrm y^+)&=-\log \pb(\mathrm y^+,\theta|\mathrm y^-,\eta,\sigma^2)\nn\\ &=
    \frac{N}{2\sigma^2}V_N(\theta)+\frac{1}{2}\norm{\theta}^2_{K^{-1}_\eta}+\frac{1}{2}\log{|K_\eta|}+\kappa
}
where $\kappa$ is a constant not depending on $\theta$ and $\eta$. Therefore, the minimization of the objective function in \eqref{Theta_hat} is equivalent to minimize (\ref{neg_lik}) and thus $\hat \theta$ also represents the Maximum A Posteriori (MAP) estimator of $\theta$.

\subsection{Hyperparameters estimation}\label{Hyperparam estim}
 In order to compute the MAP estimator $\hat \theta$, estimates of $\eta$ and $\sigma^2$ are needed. An estimate of $\sigma^2$ is obtained using a low-bias ARX model as suggested in \cite{GOODWIN_1992}.
Then, an estimate of $\eta$ is given minimizing the negative log-marginal likelihood
\begin{equation}\label{likelihood_minimization}
    \hat{\eta}=\argmin_{\eta\in \Ccal}\ell(\eta;\mathrm  y^+)
\end{equation}
where $\ell(\eta;\mathrm y^+)=-\log \pb(\mathrm y^+|\mathrm y^-,\eta,\sigma^2)$, 
\begin{equation}\label{integral}
    \pb(\mathrm y^+|{\mathrm y^-,}\eta,\sigma^2)=\int_{\Rbb^{2T}} \pb (\mathrm y^+,\theta|\mathrm y^-,\eta,\sigma^2)d\theta
\end{equation}
and $\Ccal$ is the set of constraints. For instance, if we take the TC kernels in (\ref{defTC}) then we have $\Ccal:=\{\,\eta \hbox{ s.t. } \lambda_b,\lambda_c\in(0,\infty),\; \beta_b,\beta_c\in(0,1)\}$. However, the difficulty in the optimization problem (\ref{likelihood_minimization}) is the evaluation of the objective function $\ell$ because it is difficult to find a closed form expression for the integral in \eqref{integral}. In what follows we propose a method to evaluate approximately the function $\ell(\eta;\mathrm y^+)$, given the hyperparameters vector $\eta$.
More precisely, by means of the Laplace approximation \citep[Chapter 27]{mackay2003information}, it is possible to show that an approximation of \eqref{integral} is given by 
\begin{equation}\label{lhat}
\begin{split}
\hat{\ell}(\eta ; \mathrm y^+ ) :
=&\frac{N}{2\sigma^2}V_{N}(\tilde{\theta}_{\eta})+\frac{1}{2}\norm{\tilde{\theta}_{\eta}}^2_{K_\eta^{-1}}+\frac{1}{2}\log|K_\eta|\\
&+\frac{1}{2}\log\Bigr|\frac{\partial^2{h(\theta, \eta)}}{\partial{\theta^2}} \Bigr|_{\theta=\tilde{\theta}_{\eta}}\Bigr|+\check{\kappa}\end{split}
\end{equation}
where \begin{equation} \label{eq:h}
    h(\theta, \eta):=\frac{1}{2\sigma^2}V_N(\theta)+\frac{1}{2N}\norm{\theta}^2_{K^{-1}_\eta}+\frac{1}{2N}\log{|K_\eta|},
\end{equation} $\tilde{\theta}_\eta$ is its point of minimum
and $\check \kappa$ is a constant not depending on $\eta$. 
The derivation of the aforementioned approximation will be left for a forthcoming publication.

Therefore, we can approximate  $$\ell(\eta;\mathrm y^+)\approx \hat{\ell}(\eta; \mathrm y^+).$$ \\

It is worth noting that $\mathrm y^-$ is never completely known. 
A solution to handle the initial conditions consists of setting its unknown components to zero. Such choice introduces an error which goes to zero as $N$ increases, see \cite[Section 3.2]{LJUNG_SYS_ID_1999}. For sake of completeness we report below the algorithm implementing the proposed procedure for the evaluation of $\hat{\ell}(\eta ; \mathrm y^+)$ given $\eta$.

    \begin{algorithm}
    \caption{Evaluation of $  \hat{\ell}(\eta ; \mathrm y^+ )$}\label{Alg}
     \textbf{Input}: $\cal D$, $\eta$, $\sigma^2$.\\
  \textbf{Output}: $\hat{\ell}(\eta ; \mathrm y^+ ,\tilde \theta_{\eta})$
  
    \begin{algorithmic}[1]
      \State Compute $\tilde{\theta}_\eta= \argmin_{\theta \in \Theta}V_N(\theta)+ \frac{\sigma^2}{N}\norm{\theta}^2_{{K_{\eta}^{-1}}}$
        \State Compute $\hat{\ell}(\eta ; \mathrm y^+)$ in (\ref{lhat})
    \end{algorithmic}
    \end{algorithm}

Therefore, the estimator of ${\eta}$ can be obtained by minimizing the Laplace approximation of $\ell$:
\begin{equation}\label{hat_eta}
     {\hat{\eta}}=\argmin_{\eta\in \Ccal} \hat{\ell}(\eta ; \mathrm y^+).
\end{equation}
The entire estimation procedure to obtained the model is summarized in Algorithm \ref{algo2}. It is worth noting that two optimization problems are involved in the procedure, i.e. (\ref{hat_eta}) and (\ref{Theta_hat}). The former is solved numerically using the \texttt{fmincon} Matlab routine, while the latter the \texttt{armax} routine of the Matlab System Identification Toolbox.
    \begin{algorithm}
    \caption{Estimation procedure}\label{algo2}
     \textbf{Input}: $\cal D$\\
  \textbf{Output}: $\hat \theta$
  
    \begin{algorithmic}[1]
      \State Compute an estimate $\hat \sigma^2$ of $\sigma^2$ using a low-bias ARX model 
        \State Compute $\hat \eta$ solving (\ref{hat_eta}) with $\sigma^2=\hat \sigma^2$
        \State Compute $\hat \theta$ solving (\ref{Theta_hat}) with $\eta=\hat \eta$ and $\sigma^2=\hat \sigma^2$
    \end{algorithmic}
    \end{algorithm}

\section{Numerical experiments} \label{sec num exp}
To test the proposed estimation procedure 
we performed a Monte Carlo simulation by generating 200 SISO 40-th order discrete time random models through the \texttt{drss} Matlab routine and such that the absolute value of their dominant poles belongs to the interval $\begin{bmatrix} 0.8 & 0.9\end{bmatrix}$. We fed each model with a square wave input signal whose period is equal to $650$ samples, with $50\%$ duty cycle and range between $[-0.5,0.5]$. 
We collected a total of $N=2000$ samples for the output and input signals to form the training dataset $\Dcal:=\{(y(t),u(t)),\; t=1\ldots N\}$. Subsequently, 
the models were fed with white Gausian noise with unit variance and   $N_v=1000$ samples for each model were collected to form  the validation dataset $\Dcal_v:=\{(y_v(t),u_v(t)),\; t=1\ldots  N_v\}$. The noise processes in both the experiments were chosen such that the signal to noise ratio is equal to 1.


Finally, for each training dataset we estimate the model using different estimators. It is worth noting that the input used for the training is chosen weakly exciting in order to make the identification task more challenging.
The practical length $T$ of the impulse responses is set equal to $50$.

We test the performance of each estimator using two different indices: the Coefficient Of Determination (COD) and the Average Impulse Response fit (AIR). The COD index is defined as:
\begin{equation}
 COD=100\left(1-\frac{||y_v-\hat{y}_v||^2}{||y_v-\bar y_v||^2}\right),
\end{equation}
where $y_v=[\, y_v(1) \; y_v(2)\ldots y_v(N)\,]^\top$, $\hat{y}_v$ is the vector containing the one step ahead predictions of $y_v$ using the estimated model and \alg{\bar y_v=\frac{1}{N}\sum_{k=1}^N y_v(t).} 
Let $\mathrm b=\begin{bmatrix}
    b_0& \dots& b_{T-1}
\end{bmatrix}^\top$ be the vector containing the first $T$ coefficients (i.e. the ones that are not negligible) of the actual impulse response $\{ b_k,\;\  k\in\Nbb_0\}$, and $\hat{b}=\begin{bmatrix}
    \hat{b}_0 & \dots&\hat{b}_{T-1}
\end{bmatrix}^\top$ be the vector of the estimated ones. The average impulse response fit for $\hat b$ is defined as
\begin{equation}
AIR(b,\hat{b})=100\left(1-{\frac{\sum_{k=0}^{T-1}\left(b_k-\hat{b}_k\right)^2}{\sum_{k=0}^{T-1}\left(b_k-\bar{b}\right)^2}}\right),
\end{equation}
where $$\bar{b}=\frac{1}{T}\sum_{k=0}^{T-1}b_k.$$ The average impulse response fit for $\hat c$, say $AIR(c,\hat{c})$, is defined in a similar way. Accordingly, the average impulse response fit of the estimator is 
\alg{\label{air}AIR=\frac{1}{2}(AIR(b,\hat{b})+AIR(c,\hat{c})).}

In the Monte Carlo simulation we consider the following estimators: 
\begin{itemize}
\item \textbf{TC} which denotes the estimator proposed in Section~\ref{Proposed Estimator section} equipped with the TC kernels in (\ref{defTC});
\item \textbf{D2} which denotes the estimator proposed in Section~\ref{Proposed Estimator section} equipped with the DC2 kernels in (\ref{dDC2});
\item \textbf{PEM+or} which denotes classical PEM approach to estimate Box-Jenkins (BJ) models, as implemented in the \texttt{bj.m} function of the MATLAB System Identification Toolbox, and of the form 
\begin{equation}
    y(t)=\frac{M(z)}{N(z)}u(t-1)+\frac{P(z)}{Q(z)}e(t)
\end{equation}
 with
\alg{M(z)=\sum_{k=0}^{k_1-1} m_k z^{-k},\quad N(z)=1+\sum_{k=1}^{k_1} n_kz^{-k} }\\
\alg{P(z)=1+\sum_{k=1}^{k_2} p_k z^{-k},\quad Q(z)=1+\sum_{k=1}^{k_2} q_kz^{-k} }

where the orders $k_1,k_2\in\{1,2,\dots, 10\}$. This estimator is equipped with an oracle for the fit of the average impulse response. To be specific, for any dataset, this ideal tuning provides an upper bound for AIR by selecting those model orders $k_1$ and $k_2$ 
that maximize AIR in (\ref{air});
\item \textbf{PEM+BIC} denotes the estimator defined as above, but the model orders $k_1$ and $k_2$ are chosen according to the BIC criterion, see \cite{SCHWARZ_1978}.
\end{itemize}

\begin{figure}[h]
\centering
\includegraphics[width=0.45\textwidth]{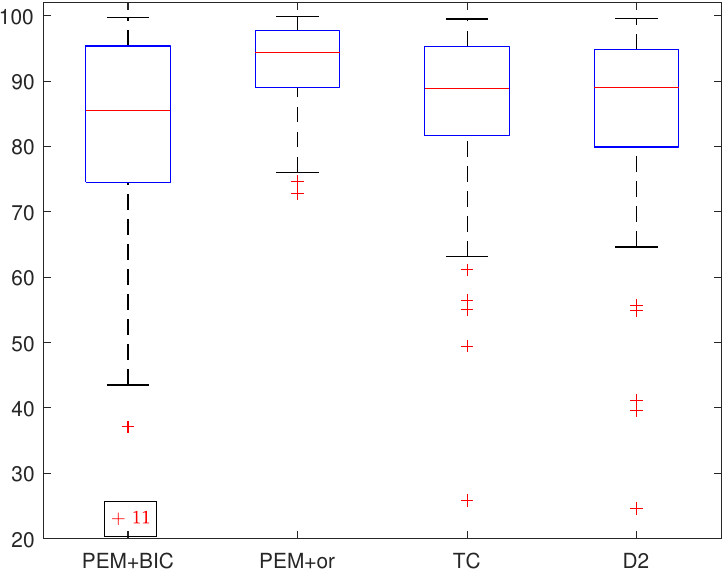}
\caption{\label{fig:air image} Boxplots of the AIR index over 200 simulations. The number in the box denotes the number of outliers not shown in the picture.}
\end{figure}

Figure \ref{fig:air image} shows the boxplots of the AIR index for the considered estimators. 
\begin{figure}[h]
\centering
\includegraphics[width=0.45\textwidth]{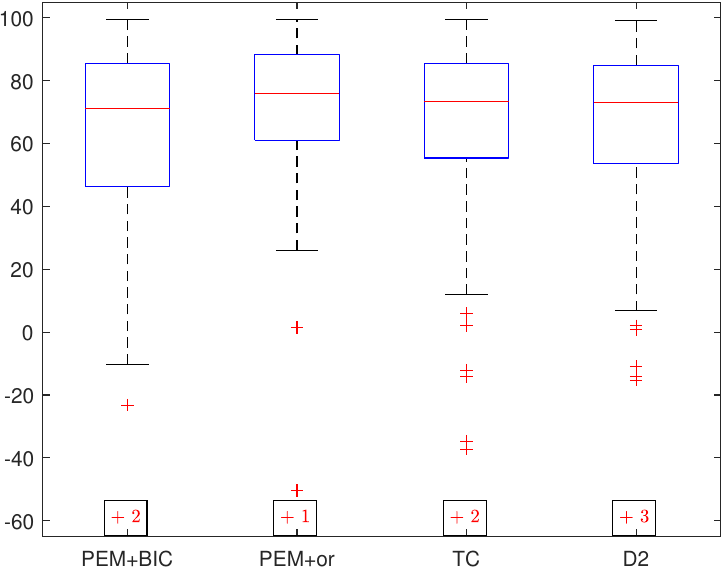}
\caption{\label{fig:cod image} Boxplots of the COD index over 200 simulations. The number in the box denotes the number of outliers not shown in the picture.}
\end{figure}
The proposed estimators, which are \textbf{TC} and \textbf{D2}, outperform \textbf{PEM+BIC}. In particular, it is possible to notice that the lower whisker of the boxplot relative to \textbf{PEM+BIC}, which indicates the lowest fit in the experiments excluding any outlier, is around 45 which is much lower with respect to the ones of \textbf{TC} and \textbf{D2} estimators (around 65). Furthermore, it is worth noting the presence of several outliers under the 20 fit percentage for the \textbf{PEM+BIC} case. As expected the ideal estimator \textbf{BIC+or} gives the best result because it selects the model knowing the true impulse response (i.e. an information that in practice one has not got).\\Figure~\ref{fig:cod image} shows the boxplots of the COD index for the considered estimators. \textbf{TC} and \textbf{D2} exhibit a better performance than the one of \textbf{PEM+BIC}. We conclude that the proposed estimator represents a valid alternative to the classic PEM estimator.

\section{Conclusions} \label{conclusion}
In this paper we have presented a novel kernel-based PEM method for the identification of a class of high-order MAX 
models which approximate a class of nonparametric forward models. The most challenging task is the estimation of the kernel hyperparameters through the minimization of the negative log-marginal likelihood because the latter lacks of a closed expression. We have proposed a procedure to evaluate the negative log-marginal likelihood by means of the Laplace approximation. The numerical results outlined in the previous section highlight the goodness of the proposed method.

There remains still open some points that we will address in future works.
More precisely, the considered high-order MAX
model class approximates a class of nonparametric models. Although the proposed estimator represents the practical way to estimate such nonparametric models, one wonder whether problem (\ref{Theta_hat}) is well posed in the infinite dimensional case, i.e. when $T=\infty$. In this case $\theta$ is a $\ell_2$ sequence and the existence of a solution to (\ref{Theta_hat}) needs further investigations. Finally, the estimation procedure has been coded with the aim of testing the performance of the estimator in the numerical experiments described in the previous section, but the efficiency in the use of computational resources (memory and time) has not been taken into account. We believe the evaluation of $\hat \ell$ can be performed efficiently exploiting some ``symmetries'' characterizing the problem.

\end{document}